\documentclass[english,11pt]{article}
\usepackage{graphicx} 
\usepackage{tikz}
\usepackage{tikz-3dplot}
\usepackage{pgfplots}
\usepackage{geometry}
\usepackage{multirow}
\usepackage{diagbox}
\usepackage{amsmath, bm}
\usepackage{amssymb}
\usepackage{amsthm}
\usepackage{booktabs}
\usepackage{varwidth}
\usepackage{array}
\usepackage{tabularray}
\PassOptionsToPackage{normalem}{ulem}
\usepackage{ulem}
\usepackage{xurl}
\makeatletter
\usepackage{xcolor}
\usepackage{color}
\usepackage{algorithm}
\usepackage[noend]{algpseudocode}
\usepackage{graphicx}
\usepackage{upgreek}
\usepackage{algpseudocode}
\newlength\fwidth
\makeatletter
\renewcommand{\ALG@beginalgorithmic}{\small}

\algrenewcommand\textproc{} 
\makeatother

\theoremstyle{remark}

\newtheorem*{rem*}{Remark}

\providecommand{\keywords}[1]
{
  \small	
  \textbf{\textit{Keywords---}} #1
}

\providecommand{\ams}[1]
{
  \small	
  \textbf{\textit{AMS subject classifications---}} #1
}

\newcommand{\cp}{\mathrm{cp}}
\newcommand{\cpbar}{\bar{\mathrm{cp}}}
\renewcommand{\S}{\mathcal{S}}

\DeclareMathOperator*{\argmin}{arg\,min}

\newtheorem{principle}{Principle}

\graphicspath{{figures/}}

\numberwithin{equation}{section}

\let\oldnl\nl
\newcommand{\nonl}{\renewcommand{\nl}{\let\nl\oldnl}}
\makeatother

\begin{document}
\title{The Closest Point Heat Method for Solving Eikonal Equations on Implicit Surfaces}
\author{Tony Wong\thanks{Department of Mathematics, University of
California, Los Angeles, USA} \and
Shingyu Leung\thanks{Department of Mathematics, The Hong Kong University of Science and Technology, Clear Water Bay, Hong Kong} \and
Byungjoon Lee\thanks{Department of Mathematics, The Catholic University of Korea, Republic of Korea}\textsuperscript{\ \ ,\P }}
\maketitle
\begingroup
\renewcommand\thefootnote{\fnsymbol{footnote}}
\footnotetext[5]{Corresponding author. Email: blee@catholic.ac.kr}
\endgroup
\maketitle

\begin{abstract}
We introduce the Closest Point Heat Method (CPHM), a novel approach for solving the surface Eikonal equation on general smooth surfaces. Building on the strengths of the classical heat method, such as simplicity of implementation and computational efficiency, CPHM integrates closest point techniques to reduce dependence on surface meshes. This embedding framework naturally extends the heat method to implicit surfaces while preserving both its efficiency and intrinsic geometric properties. Numerical experiments on benchmark geometries confirm the accuracy and convergence of the proposed method and demonstrate its effectiveness on complex shapes.
\end{abstract}

\keywords{Closest point method, Heat method, Surface Eikonal equations}

\ams{58J05, 65M06, 65M20, 65N06, 65N40, 65D18}

\section{Introduction}
The Eikonal equation is a fundamental nonlinear partial differential equation that arises in a variety of contexts involving wavefront propagation, most notably in geometric optics and computational geometry. On a smooth surface $\S$ embedded in $\mathbb{R}^3$, this equation takes the form of the surface Eikonal equation, which underpins the analysis of high-frequency surface wave propagation \cite{grimshaw1968propagation}, and is given by:
\begin{align}\label{eqn:surf_eikonal}
\begin{split}
\|\nabla_{\S} \phi(\bm{y}) \| &= \frac{1}{F(\bm{y})}, \quad \bm{y} \in \S \setminus P, \\
\phi(\bm{y}_0) &= 0, \quad \bm{y}_0 \in P,
\end{split}
\end{align}
where $\|\cdot\|$ is the Euclidean norm, $P$ is the set of source points, $F(\bm{y})$ is the wave speed, and $\phi(\bm{y})$ denotes the shortest traveling time of the wave from the source set $P$ to a point $\bm{y} \in \S$. The surface gradient $\nabla_\S \phi$ is defined by
\[
\nabla_\S \phi(\bm{y}) = \nabla \phi -\left(\textbf{N}\cdot\nabla \phi\right)\textbf{N}
\]
where $\mathbf{N}$ is the unit normal vector to the surface $\S$. Beyond wave propagation, the surface Eikonal equation is a fundamental tool in computer graphics, geometric processing, and machine learning applications such as surface reparameterization \cite{memoli2001fast}, image segmentation on manifolds \cite{liu2013splitting}, and intrinsic distance-based feature extraction \cite{bronstein2017geometric}. While the general Eikonal equation allows for spatially varying wave speed, the present work focuses on the special case when $F\left(\bm{y}\right)=1$, corresponding to the unit-speed formulation $\|\nabla_\S\phi\|=1$. This simplification is sufficient for computing geodesic distance functions and enables the development of efficient projection-based solvers.

Despite its wide applicability, computing numerical solutions to the surface Eikonal equation on general surfaces remains challenging due to the need to respect the underlying geometry and the nonlinearity of the PDE. Following the work of \cite{kimmel1998computing}, approximation-based approaches such as variants of fast marching \cite{tsitsiklis2002efficient} and fast sweeping \cite{zhao2005fast} methods have been widely adapted to surface domains. For instance, a triangulation-invariant method for anisotropic geodesic map computation improves robustness to mesh irregularities and anisotropy \cite{yoo2012triangulation}. Fast sweeping techniques have also been extended to triangulated surfaces for geodesic computation \cite{xu2008fast}, and even to implicit surfaces \cite{wong2016fast}. In addition, a number of methods based on parametric surface representations have been proposed to enhance computational efficiency: these include weighted distance map computation \cite{bronstein2007weighted}, efficient Eikonal solvers on parametric manifolds \cite{spira2004efficient}, and parallel algorithms for approximating distance maps \cite{weber2008parallel}. An optimal control approach \cite{huynh2024scalable} has also been developed, employing the Hopf–Lax formula to recast the Eikonal equation \cite{chow2019algorithm,lee2017revisiting} as a variational problem solved via convex optimization techniques. Together, these methods form a diverse set of tools for surface Eikonal solvers, though many rely on surface discretization quality and parametrization, which can limit robustness in general geometric settings. 

An alternative approach, known as the heat method, sidesteps these difficulties by exploiting the short-time behavior of the heat equation to approximate geodesic distances. First introduced by \cite{crane2013geodesics}, the heat method solves a pair of linear problems,  heat diffusion followed by the solution of a Poisson equation, making it efficient, robust to mesh quality, and naturally compatible with intrinsic surface geometry. However, the heat method typically relies on a mesh or parametrization of the surface, which can limit its applicability to more general implicit surfaces or point cloud data where such discretizations are unavailable or unreliable. 

In this paper, we devote ourselves to extending the heat method through the use of closest point techniques, which allow for solving surface partial differential equations (PDEs) without requiring an explicit parametrization or triangulated mesh. The closest point method (CPM) embeds the surface problem into a higher-dimensional Cartesian space by extending functions off the surface via a closest point extension, where each point in a tubular neighborhood is mapped to its nearest point on the surface. This enables the use of standard Cartesian finite difference schemes to approximate differential operators, such as gradients and Laplacians, directly on the surface. By avoiding mesh generation and surface fitting, this approach naturally handles surfaces represented implicitly, such as level sets or point clouds, and improves robustness in complex geometric settings. Through this embedding framework, we generalize the heat method to implicit surfaces while preserving its efficiency and intrinsic geometric fidelity.

We acknowledge that similar ideas have been explored in prior work. In particular, the authors of~\cite{king2024closest} proposed a CPM-based variant of the heat method that incorporates internal boundary conditions (IBCs) to enforce Dirichlet constraints near the source point. Their formulation modifies the right-hand side of the heat equation using a smoothed Heaviside function, requiring careful treatment of boundary interpolation and stabilization. In contrast, our approach maintains the original two-step structure of the heat method, consisting of heat diffusion followed by Poisson recovery, while leveraging explicit closest point extensions for both scalar and vector fields. This avoids the complexity of IBC enforcement and results in a simpler, mesh-free framework that naturally accommodates implicit surface representations.

\begin{figure}[t]
    \centering    
    \includegraphics[width=0.8\linewidth]{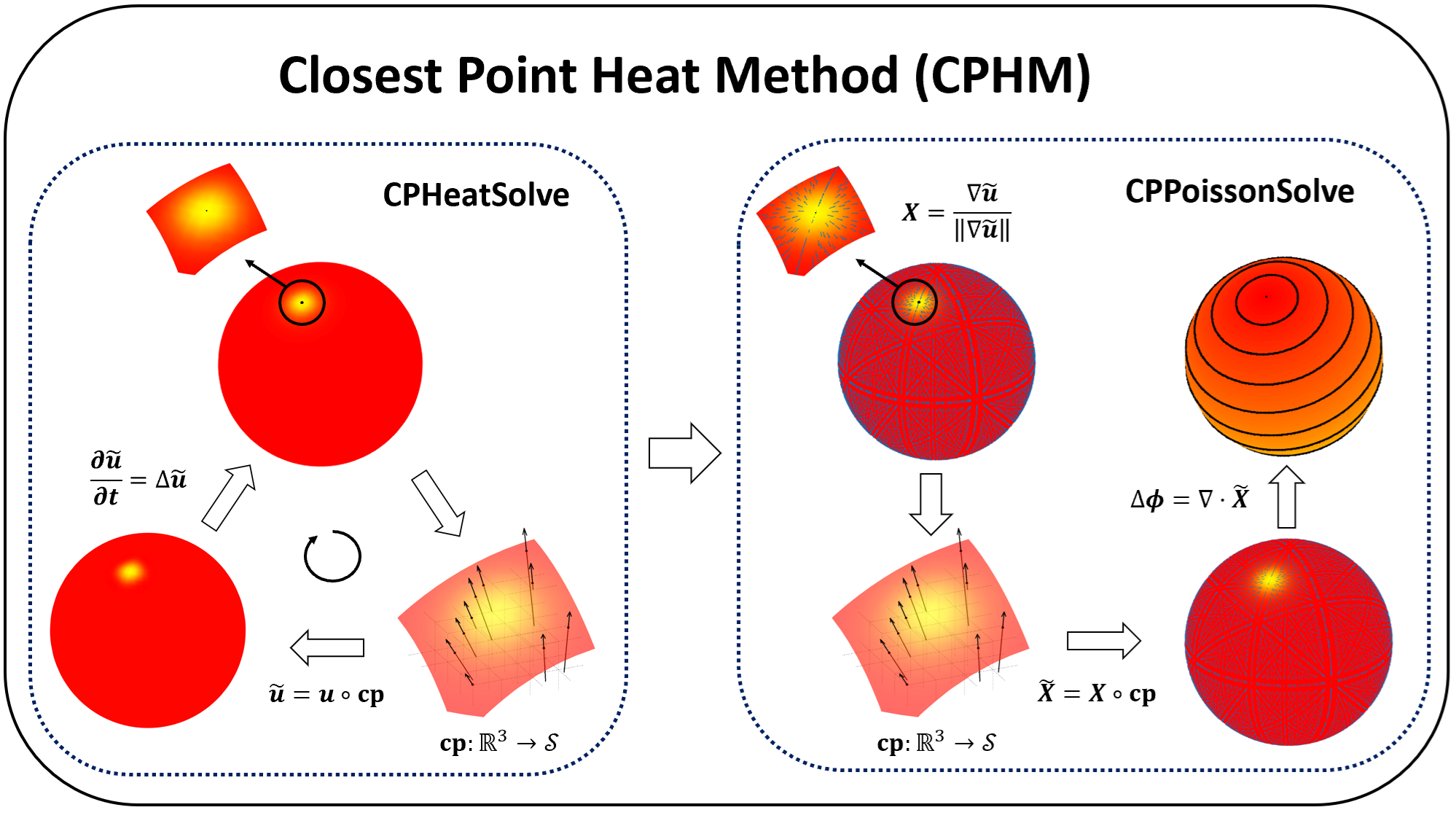}
    \caption{Overview of the proposed method. As in the original heat method~\cite{crane2013geodesics}, the approach consists of two steps: (CPHeatSolve, left) the heat equation is solved in the embedding space using the closest point extension $\tilde{u} = u \circ \mathrm{cp}$, allowing heat flow on the surface $\S \subset \mathbb{R}^3$ via $\partial \tilde{u} / \partial t = \Delta \tilde{u}$;(CPPoissonSolve, right) the normalized gradient $X = \nabla \tilde{u} / \| \nabla \tilde{u} \|$ is computed and extended, and the Poisson equation $\Delta \phi = \nabla \cdot \tilde{X}$ is then solved to recover approximate geodesic distances. All computations are carried out in $\mathbb{R}^3$ using the closest point map $\mathrm{cp}:\mathbb{R}^3 \rightarrow \S$.}
    \label{fig:overview_CPHM}
\end{figure}

Our main contributions are summarized as follows:
\begin{itemize}
\item We introduce a fast and reliable method for solving the surface Eikonal equation, called the Closest Point Heat Method (CPHM), for computing geodesic distances.
\item The proposed method eliminates the need for surface parameterizations or triangulated meshes by leveraging standard Cartesian finite difference schemes.
\item By embedding the problem in the ambient space, the method extends the classical heat method to implicit surfaces while preserving intrinsic geometric properties.
\end{itemize}
The overall procedure of CPHM is illustrated in Figure \ref{fig:overview_CPHM}.

This paper is organized as follows. Section~\ref{sec:Preliminaries} provides brief reviews of the CPM and the heat method. In Section~\ref{sec:Closest_point_heat_method}, we present the proposed method for solving the surface Eikonal equation, the Closest Point Heat Method (CPHM). Section~\ref{sec:Numerical_results} presents various numerical results that validate the performance of the method, and the final section offers conclusions and discusses directions for future work.

\section{Preliminaries}\label{sec:Preliminaries}
This section offers brief overviews of the two principal ingredients of this work: the CPM \cite{ruuth2008simple} and the heat method \cite{crane2013geodesics}.

\subsection{The Closest Point Method}
The \textit{closest point method} (CPM) is a numerical technique for solving partial differential equations (PDEs) defined on surfaces embedded in higher-dimensional spaces. The main idea of CPM is to extend the surface PDE to a neighborhood in the embedding space, enabling the use of standard Cartesian finite difference or finite element methods. We shall first introduce some key terminologies. Give a smooth surface $\mathcal{S}\subset\mathbb{R}^n$, there is a tubular neighborhood $B(\S)$ such that each point $\bm{x}\in B(\S)$ has a unique nearest point on $\S$ (in Euclidean distance). As such, we introduce the \textit{closest point function}:
\begin{equation}\label{eqn:CP_function}
\cp(\bm{x}) = \argmin_{\bm{y} \in \S} \|\bm{x} - \bm{y}\| \,, \qquad \bm{x} \in B(\S) \,,
\end{equation}
Using this function, a surface function $u: S \to \mathbb{R}$ can be extended to $B(\S)$ via the \textit{closest point extension},
\begin{equation}\label{eqn:CP_extension}
E u (\bm{x}) := u\left(\cp(\bm{x})\right) \,, \qquad \bm{x} \in B(\S)
\end{equation}
We emphasize that $Eu$ is a function that is defined in the tubular neighborhood $B(\S)$. Such extension propagates the surface function constantly along the normal direction of $\S$, and therefore allows standard differential operators (e.g., gradient, Laplacian) to act on $E u$ in the embedding space while still capturing the surface behavior. We shall state the key principles \cite{ruuth2008simple} that allow us to extend a surface PDE into the embedding space:

\begin{principle}[Gradient principle]
Let $\bm{y} \in \S$. Then 
\begin{equation}\label{eqn:gradient_principle}
\nabla_\S u\left(\bm{y}\right) = \nabla \left[ E u \left( \bm{y} \right) \right] \,.
\end{equation}
\end{principle}
\begin{principle}[Divergence principle]
Denote the surface divergence operator by $\nabla_\S \cdot$. Let $\S_\delta$ be the surface that is constantly displaced by a displacement of $\delta$ from $\S$. Suppose $\bm{v}$ is a vector field on $\mathbb{R}^3$ that is tangent at $\S$ and all $\S_\delta \subset B(\S)$. Then for $\bm{y} \in \S$, we have 
\begin{equation}
\nabla_\S \cdot \bm{v}(\bm{y}) = \nabla \cdot \bm{v}(\bm{y}) \,.
\end{equation}
\end{principle}
By applying divergence principle on $\bm{v} = \nabla_\S u = \nabla u |_{\S}$, we derive the Laplacian principle,
\begin{equation}\label{eqn:CP_laplacian}
\Delta_\S u(\bm{y}) = \Delta \left[ Eu \right](\bm{y})  \,, \qquad  \bm{y} \in \S\,.
\end{equation}
This allows us to approximate the surface Laplace–Beltrami operator $\Delta_S u$ by the standard Laplacian to $E u$ near the surface. This approach eliminates the need for surface parameterizations or triangulations and is particularly effective when combined with implicit surface representations~\cite{osher1988fronts}. For further details on CPM, see~\cite{ruuth2008simple, macdonald2010implicit}.

\subsection{The Heat Method}

To compute geodesic distances on a surface, the heat method leverages the connection between heat diffusion and shortest paths, formalized by Varadhan's asymptotic formula \cite{varadhan1967behavior}:
\begin{equation}
\phi(\bm{x}, \bm{y}) = \lim_{t \rightarrow 0} \sqrt{-4t\log k_{t,\bm{x}}\left(\bm{y}\right) }, \label{eqn:varadhan}
\end{equation}
where $\phi(\bm{x}, \bm{y})$ denotes the geodesic distance between points $\bm{x}$ and $\bm{y}$, and $k_{t,\bm{x}}(\bm{y})$ is the heat kernel, which is the fundamental solution of the heat equation initiated at $\bm{x}$.

Rather than evaluating the limit directly, the heat method approximates the distance by simulating heat flow. It begins by solving the heat equation with an initial delta impulse at $\bm{x}$:
\begin{equation*}
\partial_tu = \Delta_\S u, \quad u\left(\bm{x}, 0\right) = \delta_{\bm{x}_0}\left(\bm{x}\right),
\end{equation*}
where $\Delta_\S$ is the Laplace–Beltrami operator. The solution $ u\left(\bm{y}, t\right)$ approximates the heat kernel $k_{t,\bm{x}}\left(\bm{y}\right)$ for small $t$, providing a smooth function from which distance information can be inferred.

From the asymptotic form of the heat kernel in \eqref{eqn:varadhan}, we have
\[
u(\bm{y}, t) \approx k_{t,\bm{x}}\left(\bm{y}\right) \sim \frac{1}{(4\pi t)^{n/2}} e^{- \phi(\bm{x}, \bm{y})^2 / 4t} \cdot a(\bm{x}, \bm{y}),
\]
where $a\left(\bm{x}, \bm{y}\right)$ is a smooth, positive function. Neglecting the lower-order variation introduced by $a\left(\bm{x}, \bm{y}\right)$, we take the logarithm and differentiate with respect to $\bm{y}$ to obtain
\[
\nabla_\S u \approx u\left(\bm{y}, t\right) \cdot \nabla_\S \log k_{t,\bm{x}}\left(\bm{y}\right) \sim -\frac{\phi\left(\bm{x}, \bm{y}\right)}{2t} \nabla_\S \phi \cdot u\left(\bm{y}, t\right),
\]
which implies $\nabla_\S u \propto -\nabla_\S \phi$. In other words, the direction of steepest decrease in heat closely approximates the direction of shortest paths, that is, the geodesics on the surface.

Finally, to recover an approximation of the geodesic distance, the method solves a Poisson equation whose divergence matches that of the normalized vector field $ \bm{X} = -\nabla_\S u / \| \nabla_\S u \|$:
\[
\Delta_\S \phi = \nabla_\S \cdot \bm{X}.
\]
For more information on the heat method, we refer the reader to~\cite{crane2013geodesics}.

\section{Closest Point Heat Method}\label{sec:Closest_point_heat_method}

This section introduces a new method for solving the surface Eikonal equation, the Closest Point Heat Method (CPHM). It builds upon the classical heat method introduced by \cite{crane2013geodesics}, adapting it to the CPM framework so that all computations are carried out in the ambient Euclidean space. This formulation removes the need for surface parameterizations or triangulated meshes and instead leverages standard Cartesian finite difference schemes. In this work, we restrict ourselves to surfaces embedded in $\mathbb{R}^3$, where the closest point extension and finite difference discretization can be eifficiently implemented on regular grids. 

The method begins by solving the surface heat equation using CPM for a prescribed time $T$. As suggested in \cite{crane2013geodesics}, we set $T = (\Delta x)^2$ throughout the paper. Since the geodesic distance on the surface, the solution to the surface Eikonal equation, can be approximated by solving a related Poisson equation, the second step recovers the distance function by solving this equation within the CPM framework, with careful attention to the extension and projection steps. Algorithm~\ref{alg:CPHM} provides a summary of the CPHM procedure for a single source point. The algorithm can be easily generalized to multiple source points, which will be discussed in Sec.~\ref{sec:local_recon}.

\begin{algorithm}
\caption{Closest Point Heat Method for the Surface Eikonal Equation}
\label{alg:CPHM}
\begin{algorithmic}[1]
\State \textbf{Input:} Grid points $\bm{x} \in \mathbb{R}^3$, surface $\S$, closest point map $\operatorname{cp}$, time step $\Delta t$, penalty parameter $\gamma$, interpolation orders $p,q$
\vspace{0.5em}
\State \textbf{Initialize:} $\mathsf{u}_0 \gets \delta_{\operatorname{cp}(\bm{x}_0)}$ for a source point $\bm{x}_0 \in \S$
\State Compute extension matrices $\mathsf{E}_p$, $\mathsf{E}_q$, discrete Laplacian matrix $\mathsf{L}$, and finite difference derivative matrices $\mathsf{D}_{x_1}$, $\mathsf{D}_{x_2}$, $\mathsf{D}_{x_3}$ that approximate $\partial_{x_1}$, $\partial_{x_2}$, $\partial_{x_3}$, respectively.
\State Initialize identity matrix $\mathsf{I}$
\State $\mathsf{u}_1   \gets$ \Call{\texttt{CPHeatSolve}}{$\mathsf{u}_0$, $\Delta t$, $\gamma$, $\mathsf{I}$, $\mathsf{E}_p$, $\mathsf{E}_q$, $\mathsf{L}$} \Comment{Algorithm~\ref{alg:HeatSolver}}
\State $\phi \ \gets$ \Call{\texttt{CPPoissonSolve}}{$\mathsf{u}_1$, $\gamma$, $\mathsf{I}$, $\mathsf{E}_p$, $\mathsf{E}_q$, $\mathsf{L}$, $\mathsf{D}_{x_1}$, $\mathsf{D}_{x_2}$, $\mathsf{D}_{x_3}$} \Comment{Algorithm~3}
\State \textbf{Output:} $\phi$ restricted to $S$
\end{algorithmic}
\end{algorithm}

\subsection{Closest Point Treatment of Surface Equations}

To solve PDEs on a surface $\S$ without requiring explicit surface parameterizations, we use the CPM. This approach extends a surface PDE on the embedded surface $\S \in \mathbb{R}^3$ into a narrow tubular neighborhood $B(\S)$ where the closest point function \eqref{eqn:CP_function} is well-defined. Here, we give a brief review on the CPM for a screened Poisson equation,
\begin{equation}\label{eqn:surf_spoisson}
\Delta_\S u \left(\bm{y}\right) - cu\left(\bm{y}\right) = f\left(\bm{y}\right) \,, \qquad \bm{y} \in \S \,,
\end{equation}
where $c > 0$ and $\Delta_\S$ denotes the Laplace–Beltrami operator on $\S$. Firstly, we denote the closest point extension of the surface solution as $\tilde{u}(\bm{x}) = Eu(\bm{x})$ for $\bm{x} \in B(\S)$. Similarly, $\tilde{f} = E f$ in $B(\S)$.
Since $u \equiv \tilde{u}$ and $f \equiv \tilde{f}$ on $\S$, we can replace $u$ by $\tilde{u}$ and $f$ by $\tilde{f}$ in \eqref{eqn:surf_spoisson}:
\begin{equation}\label{eqn:surf_spoisson2}
    \Delta_\S \tilde{u} \left( \bm{y} \right) - c u\left(\bm{y}\right) = \tilde{f}\left( \bm{y} \right) \,, \qquad \bm{y} \in \S \,.
\end{equation}
Note that the closest point extension is idempotent, we have $\tilde{u} = E\tilde{u}$ in $B(\S)$. Therefore, we have
\begin{equation}\label{eqn:surflap_2_lap}
    \Delta_\S \tilde{u} = \Delta_\S \left( E \tilde{u}\right)  = \Delta \tilde{u} \mid_\S \enspace \mbox{on} \enspace \S \,,
\end{equation}
where the second equality is due to the Laplacian principle \eqref{eqn:CP_laplacian}. With \eqref{eqn:surflap_2_lap}, the surface equation \eqref{eqn:surf_spoisson2} can be expressed in terms of the Cartesian Laplacian,
\begin{equation}\label{eqn:surf_spoisson3}
\Delta \tilde{u} \left( \bm{y} \right) - cu\left(\bm{y}\right) = \tilde{f}\left(\bm{y}\right) \,, \qquad \bm{y} \in \S \,.
\end{equation}
Next, we apply the closest point extension to \eqref{eqn:surf_spoisson3}. In this way, we obtain an embedding equation for $\tilde{u}$ in $B(\S)$,
\begin{subequations}
\begin{eqnarray}
&&E\Delta \tilde{u}\left(\bm{x}\right) - c\tilde{u}\left(\bm{x}\right) = \tilde{f}\left(\bm{x}\right) \,, \qquad \bm{x} \in B(\S) \label{eqn:embedding_spoisson0_pde} \\[5pt]  \quad &&\mbox{subject to} \enspace \tilde{u}\left(\bm{x}\right) = E \tilde{u}\left(\bm{x}\right) \,, \qquad \bm{x} \in B(\S)  \,. \label{eqn:embedding_spoisson0_constraint}
\end{eqnarray}
\end{subequations}
The constraint \eqref{eqn:embedding_spoisson0_constraint} effectively enforces that the extended solution $\tilde{u}$ being a closest point extension of the surface solution $u$, and hence maintains constancy of $\tilde{u}$ along the normal direction of $\S$. Following \cite{chen2015closest}, we enforce the constraint \eqref{eqn:embedding_spoisson0_constraint} through a penalty formulation that discourages deviations from the extension condition,
\begin{equation}\label{eqn:embedding_spoisson}
E\Delta u - cu  - \gamma\left( u - E u \right) = f \enspace \mbox{in} \enspace B(\S) \,,
\end{equation}
which is referred as the \textit{embedding equation} of the surface equation \eqref{eqn:surf_spoisson}. We omit the tilde on $\tilde{u}$ and $\tilde{f}$ for notational convenience. In \eqref{eqn:embedding_spoisson}, $\gamma>0$ is a parameter that represents the penalty strength. Let $n \in \mathbb{N}$ be the dimension of the embedding space. Following \cite{von2013embedded, chen2015closest}, we choose $\gamma = 2n/(\Delta x)^2$ to balance the accuracy and effectiveness to enforce the extension condition \eqref{eqn:embedding_spoisson0_constraint}. 

Now, we review a matrix formulation (Sec.~2.2 in \cite{macdonald2010implicit}) of a finite difference scheme for the embedding equation \eqref{eqn:embedding_spoisson}. Let $\Omega_{\Delta x}$ be the collection of rectangular grid points inside the narrow band $B(\S)$. Assume $|\Omega_{\Delta x}| = N$. We introduce the vector $\mathsf{u} \in \mathbb{R}^N$ with entries $\mathsf{u}_i \approx u(\bm{x}_i)$ for each grid point $\bm{x}_i \in \Omega_{\Delta x}$.  We approximate the Cartesian Laplacian by $\Delta u \approx \mathsf{L} \mathsf{u}$, where $\mathsf{L} \in \mathbb{R}^{N \times N}$ is a Laplacian matrix that contains the weights of the standard second-order central difference scheme. When handling the closest point extension $E$ of a surface quantity, we need to assign the value at each grid point $\bm{x}_i$ to be the value at the corresponding closest point $\cp\left(\bm{x}_i\right)$. Since $\cp\left(\bm{x}_i\right)$ generally does not coincide with the grid points, we obtain the closest point value through polynomial interpolation of the surrounding grid values. Suppose the $p$-th order interpolation is used to approximate the closest point values. We may collect the interpolation weights at the stencils and incorporate them into the extension matrix $\mathsf{E}_p \in \mathbb{R}^{N \times N}$. In this way, the closest point extension $Eu$ is numerically performed by $\mathsf{E}_p \mathsf{u}$.  To balance accuracy and efficiency, we use a lower-order approximation when applying the Laplacian and higher-order interpolation for the penalty term. For example, a typical discretization of \eqref{eqn:surf_spoisson} takes the form 
\begin{equation}
\left[\mathsf{E}_p \mathsf{L} - c \mathsf{I} - \gamma\left( \mathsf{I} - \mathsf{E}_q \right) \right] \mathsf{u} = \mathsf{f} \,.
\end{equation}
where $\mathsf{I}$ is the $N$-by-$N$ identity matrix. $\mathsf{f} \in \mathbb{R}^N$ is the discretization of $f$ in $\Omega_h$. $\mathsf{E}_p$ and $\mathsf{E}_q$ denote the interpolation matrices of degree-$p$ and degree-$q$, respectively. A common choice is $p=1$, $q=3$. This strategy ensures that off-surface computations remain consistent with the surface geometry while maintaining the stability and robustness of the scheme. This section leads to a unified framework applicable to both time-dependent and steady-state surface equations in CPHM, as outlined below.

\subsubsection{Time-Dependent Case: Surface Heat Equation}

In the heat method, we have to compute a short-time solution to the surface heat equation
\begin{align}\label{eqn:surf_heat}
\begin{split}
\partial_t u\left(\bm{y},t\right) &= \Delta_\S u\left(\bm{y},t\right) \,, \quad \bm{y} \in \S, \enspace t > 0 \,, \\[5pt]
\quad u\left(\bm{y}, 0\right) &= u_0\left(\bm{y}\right) \,, \quad \bm{y} \in \S \,.
\end{split}
\end{align}
The initial condition $u_0(\bm{y})$ is taken to be the Dirac Delta distribution at the source point $\bm{y}_0$. In CPHM, we adopt a regularized Dirac Delta initial condition. See the details in Sec.~\ref{sec:local_recon}. Let $u_1\left(\bm{y}\right) \approx u\left(\bm{y}, \Delta t\right)$ be the approximated short-time heat solution. Following \cite{crane2013geodesics}, we obtain $u_1$ by applying the Backward Euler method to \eqref{eqn:surf_heat} for one time step $\Delta t$. This leads to a screened Poisson equation for $u_1$,
\begin{equation}\label{eqn:BE}
\Delta_\S u_1(\bm{y}) - \frac{1}{\Delta t} u_1(\bm{y}) = -\frac{1}{\Delta t}u_0(\bm{y}) \,, \quad \bm{y} \in \S \,.
\end{equation}
Define the closest point extensions $\tilde{u}_0 = Eu_0$ and $\tilde{u}_1=Eu_1$. The corresponding embedding equation of \eqref{eqn:BE} for $\tilde{u}_1$ is
\begin{equation}\label{eqn:embedding_BE}
    E\Delta \tilde{u}_1 - \frac{1}{\Delta t} \tilde{u}_1 - \gamma(\tilde{u}_1 - E \tilde{u}_1) = -\frac{1}{\Delta t} \tilde{u}_0 \enspace \mbox{in} \enspace B(\S) \,,
\end{equation}
which we numerically solve with the CPM. A detailed description of the numerical procedure is summarized in Algorithm~\ref{alg:HeatSolver}. 

\begin{algorithm}
\caption{\texttt{CPHeatSolve}: Closest Point Method for Surface Heat Equation}
\label{alg:HeatSolver}
\begin{algorithmic}[1]
\Function{\texttt{CPHeatSolve}}{$\mathsf{u}_0$, $\Delta t$, $\gamma$, $\mathsf{I}$, $\mathsf{E}_p$, $\mathsf{E}_q$, $\mathsf{L}$}
    \State Let $\mathsf{u}_1 \in \mathbb{R}^N$ with entries $\left[\mathsf{u}_1\right]_i \approx \tilde{u}_1\left(\bm{x}_i\right)$ in \eqref{eqn:embedding_BE}. Compute $\mathsf{u}_1$ by solving the linear system:     \begin{equation}\label{eqn:embedding_BE_discretized}
    \left[\mathsf{E}_p \mathsf{L} - \frac{1}{\Delta t} \mathsf{I} - \gamma(\mathsf{I} - \mathsf{E}_q)\right] \mathsf{u}_1 = -\frac{1}{\Delta t} \mathsf{u}_0 
    \end{equation}
    
    \State \Return $\mathsf{u}_1$
\EndFunction
\end{algorithmic}
\end{algorithm}

\subsubsection{Steady-State Case: Surface Poisson Equation}\label{sec:surf_poisson}


To obtain the Eikonal solution, we have to solve the surface Poisson equation,
\begin{equation}\label{eqn:surf_poisson}
\Delta_\S \phi = \nabla_\S \cdot \bm{X} \enspace \mbox{on} \enspace \S \,, \qquad \mbox{where} \enspace \bm{X} = -\nabla_\S u_1 / \| \nabla_\S u_1 \| \,,
\end{equation}
which corresponds to \eqref{eqn:surf_spoisson} with $c=0$ and $f = \nabla_\S \cdot \bm{X}$. In order to derive the embedding equation for \eqref{eqn:surf_poisson}, we provide some details on how to perform the closest point extension of its right-hand side function $f$. By the gradient principle \eqref{eqn:gradient_principle}, we can replace the surface gradient of $u_1$ by the Cartesian gradient of $\tilde{u}_1 = E u_1$. Therefore, we have
\begin{equation}
\bm{X} = -\nabla \tilde{u}_1 / \| \nabla \tilde{u}_1 \| \enspace \mbox{on} \enspace \S \,.
\end{equation}
Since $\bm{X}$ is parallel to $\nabla_\S u_1$, hence it is tangent to $\S$. Define $\tilde{\bm{X}}$ as the closest point extension of the vector field $\bm{X}$, i.e., $\tilde{\bm{X}}(\bm{x}) = \bm{X}\left(\cp(\bm{x})\right)$. By the Divergence principle, we can replace the surface divergence of $\bm{X}$ by the Cartesian divergence of $\tilde{\bm{X}}$ when computing
\begin{equation}
f = \nabla_\S \cdot \bm{X} = \nabla \cdot \tilde{\bm{X}} \enspace \mbox{on} \enspace \S \,.
\end{equation}
Following the derivation in the previous sections, the embedding equation of $\tilde{\phi} = E\phi$ for \eqref{eqn:surf_poisson} is
\begin{equation}
E \Delta \phi - \gamma (\phi - E \phi) = \tilde{f} \enspace \mbox{in} \enspace B(\S) \,, \quad \mbox{where} \enspace \tilde{f}\left(\bm{x}\right) = [\nabla \cdot \tilde{X} ]\left(\cp(\bm{x})\right) \,.
\end{equation}
The corresponding solution method is outlined in Algorithm~\ref{alg:PoissonSolver}.


By using the closest point extension operator and a penalty-based embedding strategy, both parabolic and elliptic surface PDEs can be treated in a consistent and efficient framework. This approach allows the use of standard Cartesian discretizations within a narrow band around the surface, simplifying implementation while preserving the intrinsic geometry of the problem. 

\begin{algorithm}
\caption{\texttt{CPPoissonSolve}: Recovery of Distance via Closest Point Poisson Solve}
\label{alg:PoissonSolver}
\begin{algorithmic}[1]
\Function{\texttt{CPPoissonSolve}}{$\mathsf{u}_1$, $\gamma$, $\mathsf{I}$, $\mathsf{E}_p$, $\mathsf{E}_q$, $\mathsf{L}$, $\mathsf{D}_{x_1}$, $\mathsf{D}_{x_2}$, $\mathsf{D}_{x_3}$}
    \vspace{5pt}
    \State Compute $\nabla \tilde{u}_1 = [\partial_{x_1} \tilde{u}_1, \,  \partial_{x_2} \tilde{u}_1, \, \partial_{x_3} \tilde{u}_1]^T$ at each grid point $\bm{x}_i$ component-by-component:
    $$\partial_{x_j}\tilde{u}\left(\bm{x}_i \right) \gets \left[\mathsf{D}_{x_j} \mathsf{u}_1\right]_i \,, \qquad \mbox{for} \enspace j=1,2,3 \,.$$
    \State Compute normalized heat gradient $\bm{X}=[X_1, X_2, X_3]^T$ at each grid point $\bm{x}_i$:
    $$ \bm{X}(\bm{x}_i) = -\nabla \tilde{u}_1(\bm{x}_i) / \|\nabla \tilde{u}_1 (\bm{x}_i)\| \,. $$    
    \State Compute $\partial_{x_j} X_j$ for each grid point $\bm{x}_i$, for $j=1,2,3$:
    $$ \partial_{x_j} X_j\left(\bm{x}_i \right) \gets \left[\mathsf{D}_{x_j} \left[X_j(\bm{x}_1), X_j(\bm{x}_2),\cdots,X_j(\bm{x}_N)\right]^T \right]_i $$
    \State Compute $f = \nabla \cdot \bm{X} = \partial_{x_1} X_1 + \partial_{x_2} X_2 + \partial_{x_3} X_3$ at each grid point $\bm{x}_i$, \Comment{Cartesian divergence}
    $$ f(\bm{x}_i) \gets \partial_{x_j} X_1\left(\bm{x}_i \right) + \partial_{x_j} X_2\left(\bm{x}_i \right) + \partial_{x_j} X_3\left(\bm{x}_i \right) $$       
    \State Let $\mathsf{f} \in \mathbb{R}^N$ with entries $\mathsf{f}_i = f(\bm{x}_i)$. Solve for $\phi$ from the embedding equation:
    \begin{equation}\label{eqn:surf_poisson_disretization}
    \left[ \mathsf{E}_p \mathsf{L} - \gamma (\mathsf{I} - \mathsf{E}_q) \right] \upphi = \mathsf{E}_q \mathsf{f}.
    \end{equation}
    \State \Return $\upphi$
\EndFunction
\end{algorithmic}
\end{algorithm}

\subsection{Approximation of the Dirac Delta Function on a Surface}\label{sec:local_recon}
In the CPHM, the accuracy of the surface delta function approximation plays a crucial role in the overall quality of the solution. Since this approximation represents a concentrated source term, any error in its support width, placement, or the underlying geodesic distance can propagate through the simulation and degrade the fidelity of results such as heat flow or distance maps.

Let $\phi_\S(\bm{x})$ be the geodesic distance between the source point $\bm{x}_0$ and $\bm{x}$ on the surface $\S$. When solving the heat equation \eqref{eqn:surf_heat}, we approximate the surface Dirac Delta initial condition by
\begin{equation}\label{eqn:dirac_approx}
\delta_{\bm{x}_0, H}(\bm{x}) = \begin{cases}
\dfrac{2\pi}{(\pi^2 - 4)H^2} \left[ 1 + \cos\left( \dfrac{\pi \, \phi_\S(\bm{x})}{H} \right) \right], \quad &\mbox{if} \quad \phi_\S(\bm{x}) \leq H , \\[10pt]
0, \quad &\mbox{if} \quad \phi_\S(\bm{x}) > H ,
\end{cases}
\end{equation}
which is obtained by generalizing a radially symmetric approximation of the Dirac delta distribution on $\mathbb{R}^2$ with compact support of radius $H$ \cite{hosseini2016regularizations}, replacing the Euclidean distance with the geodesic distance $d_\S$. 

To assign such an initial condition in a narrow band of the surface $\S$, we first search for the grid points whose closest points lie within a Euclidean ball of radius $H$ centered at $x_0$. Then, for each of these closest points, we approximate their geodesic distance from $x_0$ by fitting a bivariate polynomial surface $\tilde{\S}$ near $x_0$. The geodesic distance on the original surface is then approximated by that on the fitting surface $\tilde{\S}$, which can be computed by solving a geodesic equation \cite{kasap2005numerical}.

The local reconstructed surface $\tilde{\S}$ is defined in a local coordinate system centered at $\bm{x}_0$, where the coordinate axes are aligned with the principal directions at that point: the surface normal $\bm{N}(\bm{x}_0)$ and two orthonormal tangent vectors $\bm{t}_1(\bm{x}_0)$ and $\bm{t}_2(\bm{x}_0)$ spanning the tangent plane, as described in \cite{leung2009grid}. In this coordinate frame, any point near $\bm{x}_0$ can be expressed as $(\xi, \eta, \tilde{z})$, where $(\xi, \eta)$ are coordinates in the tangent plane and $\tilde{z}$ is the height along the normal direction.

We fit a second-order polynomial surface of the form
\[
\tilde{z}(\xi, \eta) = a_0 + a_1 \xi + a_2 \eta + a_3 \xi^2 + a_4 \xi \eta + a_5 \eta^2,
\]
using least squares fitting from the local closest point data transformed into this coordinate system. This reconstructed surface $\tilde{\S}$ captures the local curvature of $\S$ near $x_0$ and enables accurate approximation of intrinsic quantities such as geodesic distance.

\begin{figure}[t]
    \centering    
    \includegraphics[width=0.8\linewidth]{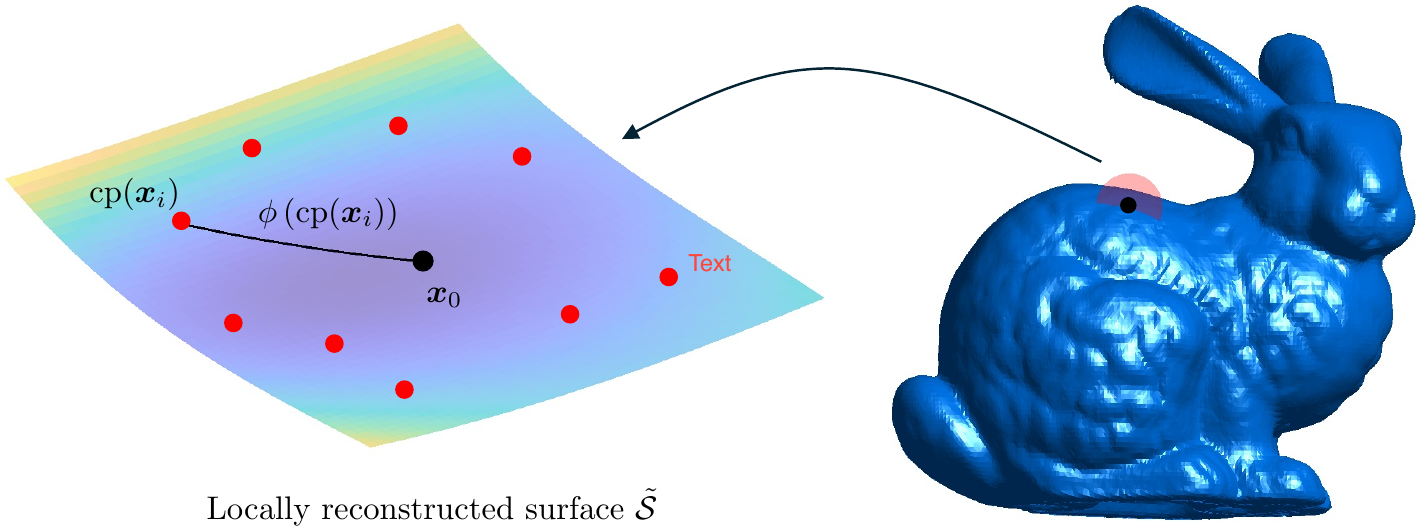}
    \caption{Here, we illustrate the local reconstruction of the surface $\S$ used to approximate the geodesic distance near the source point $\bm{x}_0$. The intersection of a Euclidean ball centered at $\bm{x}_0$ and the surface $\S$ is approximated by a bivariate polynomial surface $\tilde{\S}$ defined in a local coordinate system at $\bm{x}_0$. The geodesic distance between $\bm{x}_0$ and a closest point $\cp(\bm{x}_{i,j})$ on $\S$ is then approximated by the geodesic distance between the corresponding points on $\tilde{\S}$ within this local frame.}
    \label{fig:local_recon}
\end{figure}

Empirically, we have found that $H = 2\Delta x$ gives the optimal accuracy. Note that this surface local reconstruction procedure can be easily generalized to multiple source points, provided that their support do not overlap.

\subsection{Extension to Open Surfaces}\label{sec:open_surfaces}

Here, we consider an embedded surface $\S \subset \mathbb{R}^3$ with a smooth and co-dimension one boundary $\partial \S$. Such open surfaces pose two challenges to CPHM. Firstly, we must impose appropriate boundary conditions (BCs) for both the surface heat and Poisson equations on $\partial \S$. Secondly, we must adapt the closest point extension to the surface PDEs with the presence of BCs. Throughout this section, we let $\bm{n} = \bm{n}(\bm{y})$ be the unit outward normal at a point $\bm{y} \in \partial \S$.

\subsubsection{Correct boundary conditions}


Our objective is to compute the geodesic distance field from the source points without boundary effects; a crucial step is to compute a boundary-free heat gradient. To begin the discussion, suppose we solve the heat equation \eqref{eqn:surf_heat} with homogeneous Neumann (no-flux) BC, 
\begin{equation*}
\partial_n u = \langle \bm{n}, \nabla_\S u \rangle = 0 \enspace \mbox{on} \enspace \partial\S \,.
\end{equation*}
Under this BC, the resulting heat gradient is orthogonal to the unit outward normal of $\partial \S$. Hence, we can see that the boundary geometry alters the heat gradient, which leads to a distorted Eikonal solution that differs from the geodesic distance field. To improve this, \cite{crane2013geodesics} proposes averaging the (homogeneous) Neumann heat solution with the Dirichlet heat solution. Such heuristic modification appears to mitigate the boundary effect. Alternatively, we consider a modified Neumann BC for \eqref{eqn:surf_heat}
\begin{subequations}
\begin{equation}\label{eqn:NBC_heat_open}
\partial_n u = g \,, \enspace \mbox{where} \enspace g = \langle \bm{n}, \nabla_\S u \rangle \,.
\end{equation}
In this formalism, the heat gradient $\nabla_\S u$ is completely determined by the heat source. At the boundary, the normal derivative of the heat solution is defined by the inner product of the heat gradient and the unit outward normal of $\partial \S$. This self-consistent BC does not interfere with the heat propagation from the source points. Consequently, we obtain the heat gradient that has the boundary effects eliminated. As in Sec.~\ref{sec:surf_poisson}, we impose the Eikonal gradient to be the normalized heat gradient $\bm{X} = -\nabla_\S u/\|\nabla_\S u\|$ on $\S\cup\partial\S$. This constraint leads to the minimization problem $\min\limits_{\phi}\|\nabla_\S\phi - \bm{X}\|_2$. The Euler–Lagrange equation of this problem is the Poisson equation \eqref{eqn:surf_poisson} with the Neumann BC
\begin{equation}\label{eqn:NBC_poisson_open}
    \partial_n \phi = \langle \bm{n}, \bm{X} \rangle \enspace \mbox{on} \enspace \partial \S \,.
\end{equation}
\end{subequations}
In conclusion, we supplement the surface heat and Poisson equations with the correct BCs \eqref{eqn:NBC_heat_open} and \eqref{eqn:NBC_poisson_open}, respectively.

\subsubsection{Modified CPHM for boundary condition}

To handle an open surface, we adopt a modified closest point function \cite{macdonald2011solving} for the grid points near the boundary. Moreover, the discretization of the inhomogeneous Neumann BC requires special care such that the second-order accuracy of the Poisson solver can be retained. We leave the details in the Appendix \ref{sec:NBC}. Firstly, the discretization \eqref{eqn:embedding_BE_discretized} of the heat solution is changed to
\begin{equation}\label{eqn:embedding_BE_open}
\left(\bar{\mathsf{E}}_p \mathsf{L} - \frac{1}{\Delta t} \mathsf{I} \right) \mathsf{u}_1 - \gamma(\mathsf{u}_1 - \bar{\mathsf{E}}_q \mathsf{u}_1 - \mathsf{g})  = -\frac{1}{\Delta t} \mathsf{u}_0
\end{equation}
where $\bar{\mathsf{E}}_p$ and $\bar{\mathsf{E}}_q$ are the interpolation matrices of order $p$ and $q$ at the modified closest points, respectively. The vector $\mathsf{g} \in \mathbb{R}^N$ encodes a second-order accurate discretization of the Neumann BC \eqref{eqn:NBC_heat_open}. Since the right-hand side function $g$ of \eqref{eqn:NBC_heat_open} is a linear operator of the heat solution $u$, its corresponding discretization $\mathsf{g}$ is also a linear with respects to the discretized heat solution $\mathsf{u}_1$. Hence, we can explicitly write
\begin{equation}\label{eqn:NBC_poisson_matrix}
    \mathsf{g} = \bar{\mathsf{E}}_g \mathsf{u}_1 \,, \enspace \mbox{for some matrix} \enspace \bar{\mathsf{E}}_g \in \mathbb{R}^{N\times N} \,.
\end{equation}
Upon substituting \eqref{eqn:NBC_poisson_matrix} into \eqref{eqn:embedding_BE_open}, $\mathsf{u}_1$ now satisfies the linear system:
\begin{equation}\label{eqn:embedding_BE_open2}
\left[ \bar{\mathsf{E}}_p \mathsf{L} - \frac{1}{\Delta t} \mathsf{I}- \gamma(\mathsf{I} - \bar{\mathsf{E}}_q - \bar{\mathsf{E}}_g ) \right] \mathsf{u}_1  = -\frac{1}{\Delta t} \mathsf{u}_0
\end{equation}
Similarly, for the discrete Poisson solution $\mathsf\phi$ in \eqref{eqn:surf_poisson_disretization}, a second-order discretization of the Neumann BC \eqref{eqn:NBC_poisson_open} modifies the extension constraint to $\mathsf{\phi} = \mathsf{E}_q \mathsf{\phi} + \mathsf{g}_2$, where $\mathsf{g}_2 \in \mathbb{R}^N$ encodes a second-order accurate discretization of the Neumann BC \eqref{eqn:NBC_poisson_open}. Upon substituting the modified constraint into \eqref{eqn:surf_poisson_disretization}, we have
\begin{equation}
\left[ \bar{\mathsf{E}}_p \mathsf{L} - \gamma (\mathsf{I} - \bar{\mathsf{E}}_q) \right] \phi = \bar{\mathsf{E}}_q \mathsf{f} + \gamma \mathsf{g}_2 \,.
\end{equation}
These modifications allow CPHM to incorporate second-order accurate treatments of Neumann boundary conditions on open surfaces while maintaining the overall structure of the original algorithm. 

\begin{rem*}[On the Accuracy of CPHM]
Although CPHM incorporates higher-order accurate discretizations for interpolation and the Laplacian within the closest point framework, the overall accuracy of the method remains first-order, consistent with the original heat method \cite{crane2013geodesics}. This limitation primarily stems from the design of \texttt{CPHeatSolve} (Algorithm~\ref{alg:HeatSolver}), which adopts the same time discretization strategy, using a time step of $\Delta t = (\Delta x)^2$ to ensure stability. The resulting distance function inherits the $\mathcal{O}(\Delta x)$ error introduced by the short-time heat flow approximation, which only asymptotically recovers the gradient direction of the true geodesic distance. We also observe a potential further reduction in accuracy for certain open surface problems, primarily due to errors in the boundary condition approximation. A more detailed investigation into the accuracy of CPHM will be provided in Section~\ref{sec:convergence_study}.
\end{rem*}

\section{Numerical results}\label{sec:Numerical_results}
In this section, we present numerical results to validate the performance of the proposed method, CPHM. The first subsection is dedicated to a convergence study, which includes two parts: numerical tests on the unit sphere, a closed surface, with both single and multiple point sources; and experiments on an open surface to evaluate the method's behavior near boundaries. A series of examples on complex surfaces follows, demonstrating the versatility of CPHM in accurately capturing geodesic distances on geometries with intricate features and topology. The penalty parameter $\gamma$ is chosen to be $2n/\left(\Delta x\right)^2$, where $n$ is the dimension of the data, and this value is used consistently across all examples. All computations were implemented in MATLAB using code based on the GitHub repository at \url{https://github.com/cbm755/cp_matrices}, and executed on a personal laptop (Apple M1 Pro, 3.2 GHz processor, 16 GB memory). Every linear system arising in the CPHM algorithm (Algorithm~\ref{alg:CPHM}) is solved using MATLAB's \texttt{backslash} operator. We note that one could also adopt a multigrid approach, as in \cite{chen2015closest}, to enhance the efficiency of the solver.

\subsection{Convergence Study}\label{sec:convergence_study}

\subsubsection{Closed Surface: Unit Sphere}

\begin{figure}[t]
\centering
\includegraphics[width=0.99\textwidth]{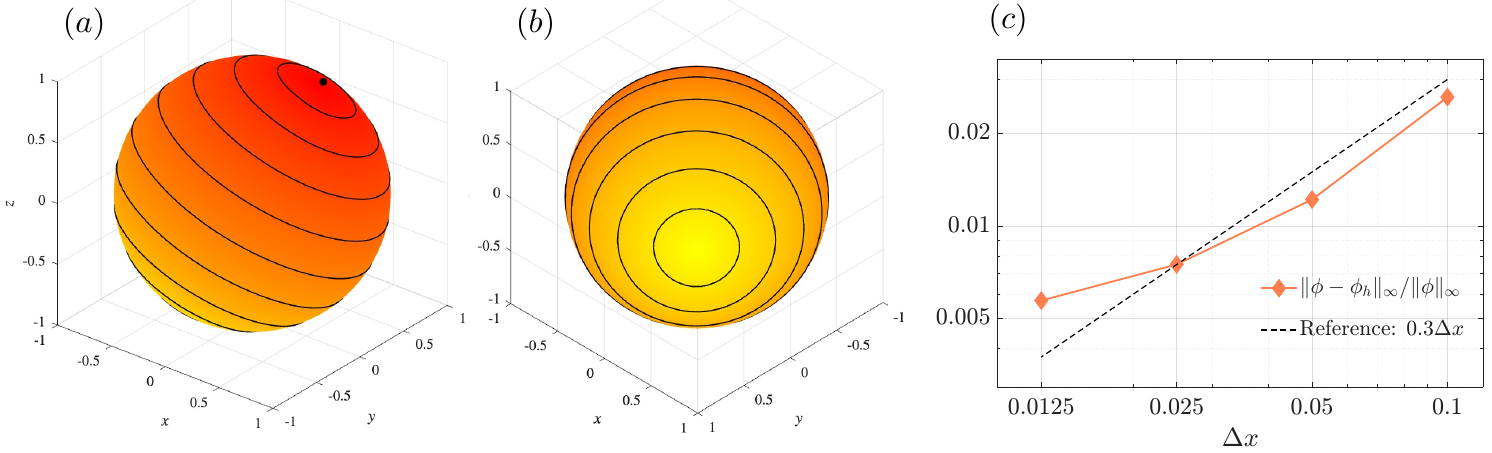}
\caption{(a) CPHM Eikonal solution obtained with $\Delta x = 0.05$ (number of grid points is 41,870, see Table \ref{tab:sphere_runtime}). The source point is located by the black marker. In (b), we show the eikonal solution from the bottom view. (c) We show the first-order convergence in the relative error $l_\infty$ norm.}
\label{fig:convg_sphere}
\end{figure}

\begin{figure}[t]
\centering
\includegraphics[width=0.99\textwidth]{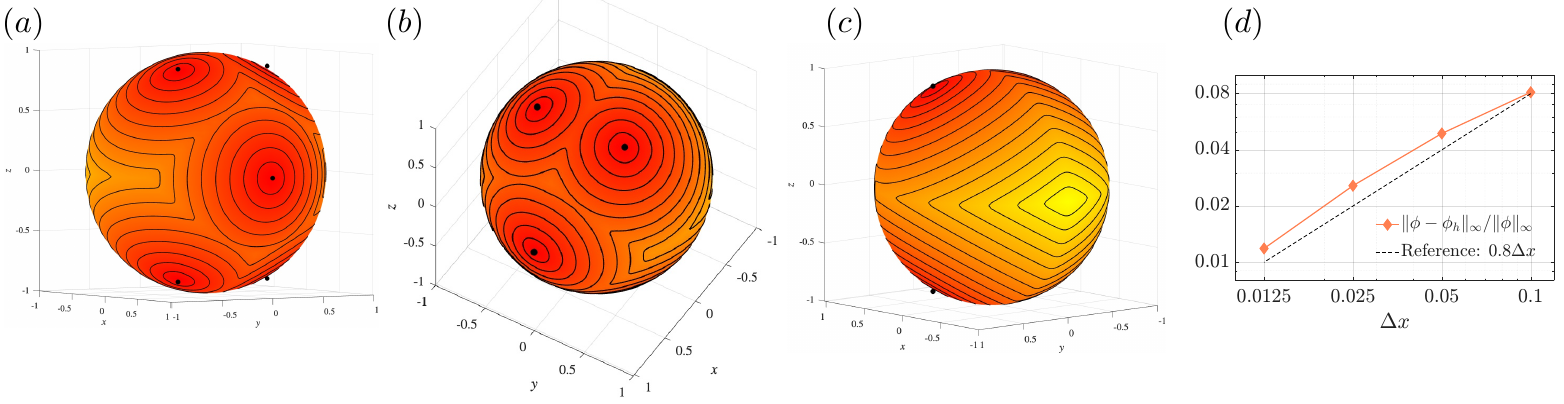}
\caption{(a), (b), (c): CPHM Eikonal solution with 5 source points, and $\Delta x = 0.025$ (number of grid points is 166,390) in different viewing angles. Each source point is located by a black marker. In (d), we show the first-order convergence in the relative errors in $l_\infty$ norm.}
\label{fig:sphere_multiple_source}
\end{figure}

We begin our convergence study of CPHM with a simple benchmark problem. Let $\S$ be the surface of the unit sphere $\| \bm{x} \| = 1$ embedded in $\mathbb{R}^3$. We adopt the spherical coordinates $(\varphi, \theta)$, where $\varphi \in [0, 2\pi)$ represents the azimuth angle and $\theta \in [0, \pi/2]$ denotes the co-latitudinal angle measured from the positive $z$-axis. In the first example, we set a single source point at $(\varphi, \theta) = (\pi/4, \pi/3)$. The Eikonal solution shown in Fig.~\ref{fig:convg_sphere}(a,b) exhibits a smooth profile and equispaced contour lines that are consistent with the exact geodesic distance on the unit sphere. The plots of relative error in Fig.~\ref{fig:convg_sphere}(c) confirm that the proposed method achieves almost first-order convergence, in agreement with the original heat method \cite{crane2013geodesics}. We give the number of grid points inside the computational band (Length($\phi^h$)) and computational times for the main steps of CPHM, including the total time for local reconstruction in approximating Dirac Delta initial condition, the heat solver, and the Poisson solver. It is noteworthy that the computational time for local reconstruction is largely independent of the mesh size $\Delta x$. It is because the support $H$ of the smoothened Dirac Delta initial condition scales with $\Delta x$ (we choose $H=2\Delta x$) In the second example, we consider five source points at $(\varphi, \theta) = (0,0), \, (\pm \, \pi/3, \, \pm \, \pi/3)$. Unlike the case of a single source point, the CPHM Eikonal solution in Fig.~\ref{fig:sphere_multiple_source}(a-c) is able to accurately capture the kinks at the intersection of characteristic curves of the Eikonal equation, supported with a first-order convergence in Fig.~\ref{fig:sphere_multiple_source}(d).

\begin{table}[htbp]
    \centering
    \caption{Computational complexity and times for Fig.~\ref{fig:convg_sphere}.}
    \resizebox{\textwidth}{!}{
    \begin{tabular}{c c c c c }
        \hline
        $\Delta x$ & Length($\phi^h)$ & Time (Local reconstruction) & Time (Heat solver) & Time (Poisson solver) \\ \hline
        0.1  & 10,906 & 0.7072 (s)   & 0.2387 (s) & 0.2193 (s) \\
        0.05  & 41,870 & 0.5982 (s)  & 1.6003 (s) & 1.6593 (s)  \\
        0.025  & 166,390 & 0.4262 (s) & 15.7687 (s)  & 16.1150 (s)  \\
        0.0125 & 663,454 & 0.4366 (s) & 306.1670 (s)  & 191.4646 (s)  \\ \hline
    \end{tabular}
    }
    \label{tab:sphere_runtime}
\end{table}

\subsubsection{Open Surfaces}

\begin{figure}[t]
\centering
\includegraphics[width=0.8\textwidth]{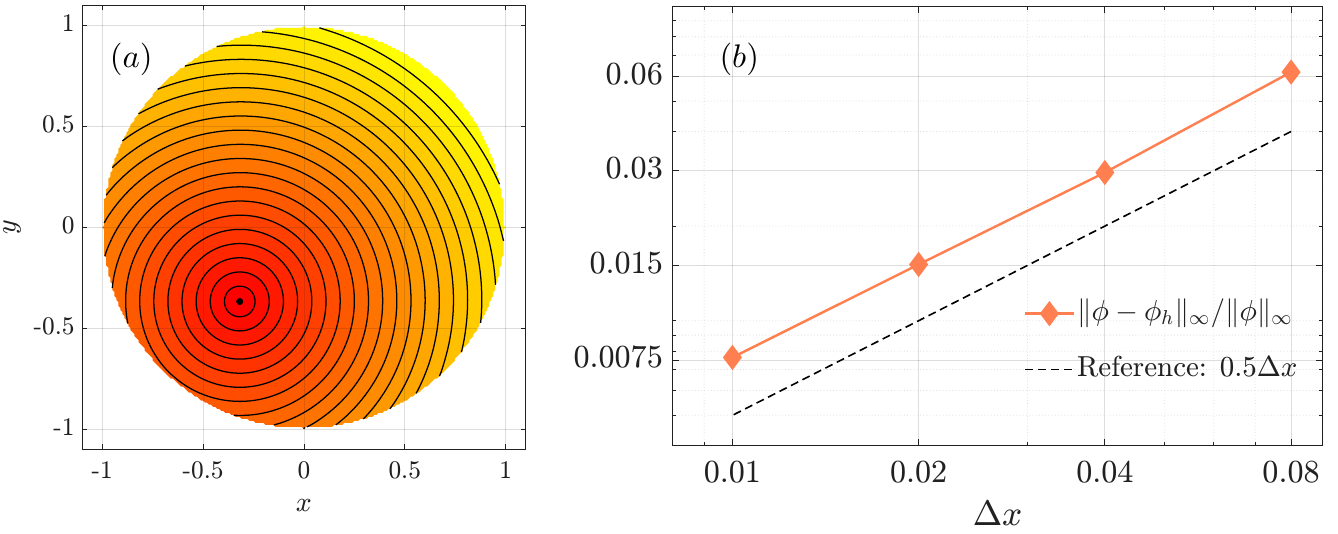}
\caption{(a): CPHM Eikonal solution of the planar unit disk with source point $\bm{x}_0=-(\pi^{-1}, e^{-1})$ obtained with $\Delta x=0.01$ (number of grid points = 33,745). (b): First-order convergence of relative errors in $L_\infty$-norm.}
\label{fig:disk}
\end{figure}

First, to validate the proposed boundary condition in Sec.\ref{sec:open_surfaces}, we apply the modified CPHM to the (planar) unit disk embedded in $\mathbb{R}^2$. We observe that the Eikonal solution in Fig.~\ref{fig:disk} propagates correctly near the boundary. The relative error in $L_\infty$-norm confirms the first-order convergence of the modified CPHM solution (Fig.\ref{fig:disk}(b)).

\begin{figure}[t]
\centering
\includegraphics[width=0.99\textwidth]{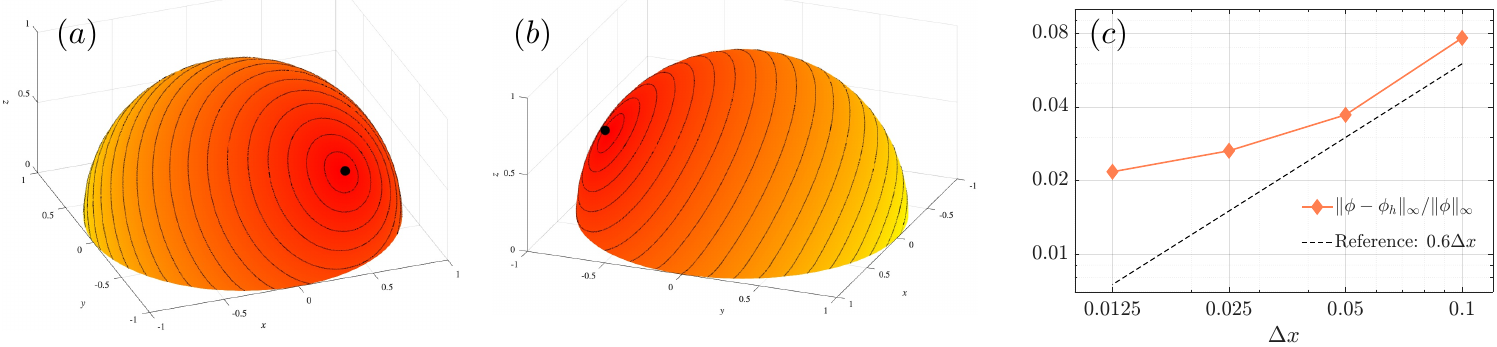}
\caption{(a), (b): CPHM Eikonal solution on the upper hemisphere with source point (shown as black marker) at $(\varphi, \theta) = (5\pi/3, \, 3\pi/10)$  obtained with $\Delta x=0.025$ (number of grid points = 89,989). (b): The relative errors in $L_\infty$-norm.}
\label{fig:hemisphere}
\end{figure}

Next, we consider the upper hemisphere $\{(x,y,z)\in\mathbb{R}^3: x^2 + y^2 + z^2 = 1, z > 0 \}$. As shown in Fig.~\ref{fig:hemisphere}(a-b), the CPHM Eikonal solution is qualitatively correct, particularly near the boundary. However, the relative error converges slower than first-order (see Fig.~\ref{fig:hemisphere}(c)). We attribute the slow convergence to error propagation from the heat solver into the Poisson solver. More precisely, the right-hand side of \eqref{eqn:surf_poisson} is computed by taking the numerical divergence of the normalized gradient of the discrete heat solution. Since the heat equation \eqref{eqn:surf_heat} is solved with a standard second-order Laplacian discretization, the resulting (normalized) gradient is only first-order accurate. Thus the Poisson right-hand side carries an $O(\Delta x)$ truncation error that limits the overall convergence rate. This accuracy bottleneck is inherent to the heat-method pipeline, not to the proposed boundary condition or the underlying closest point framework: as shown for the unit disk (Fig.~\ref{fig:disk}(b)), we do obtain first-order convergence. The boundary condition is correct at the continuous level, but a first-order convergence may require higher-accuracy discretizations in both the heat and Poisson solvers. Open surfaces therefore remain an important challenge for further improving the heat method.

\subsection{Miscellaneous Examples}

\begin{figure}[!ht]
\centering
\includegraphics[width=0.9\textwidth]{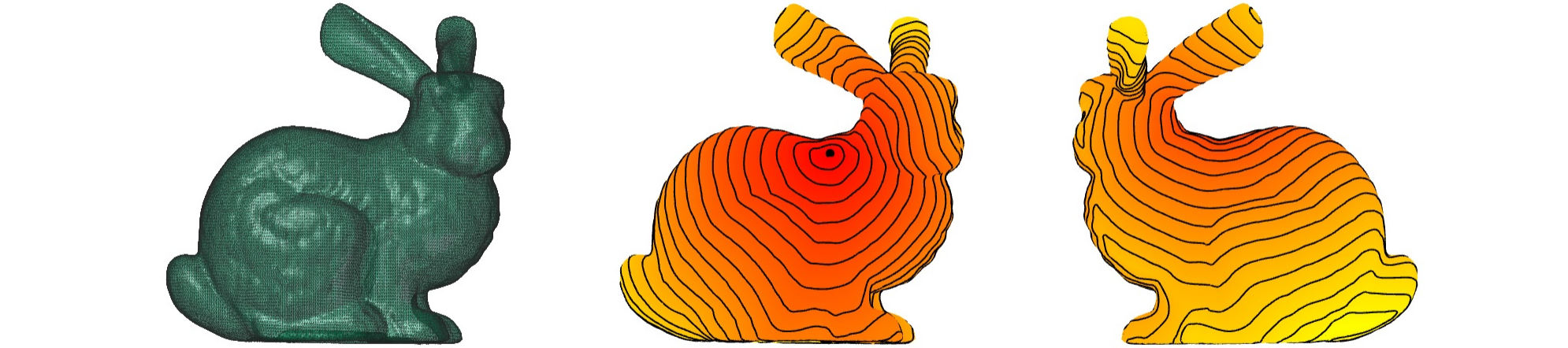}\\[10pt]
\includegraphics[width=0.9\textwidth]{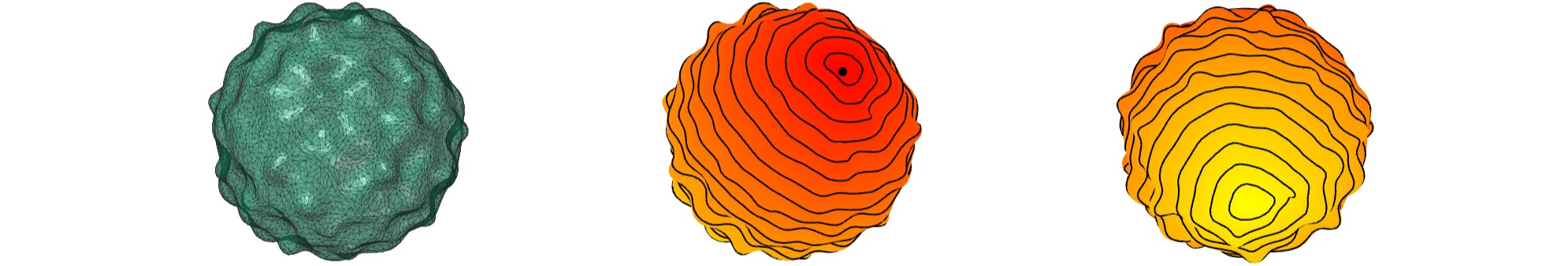}\\ [10pt]
\includegraphics[width=0.9\textwidth]{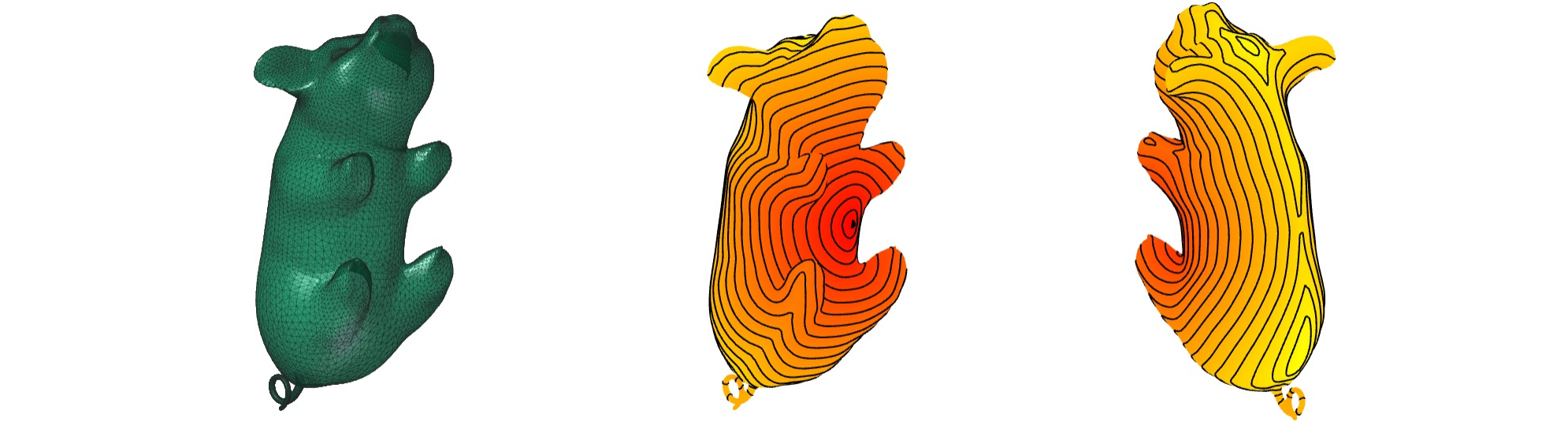} \\[10pt]
\includegraphics[width=0.9\textwidth]{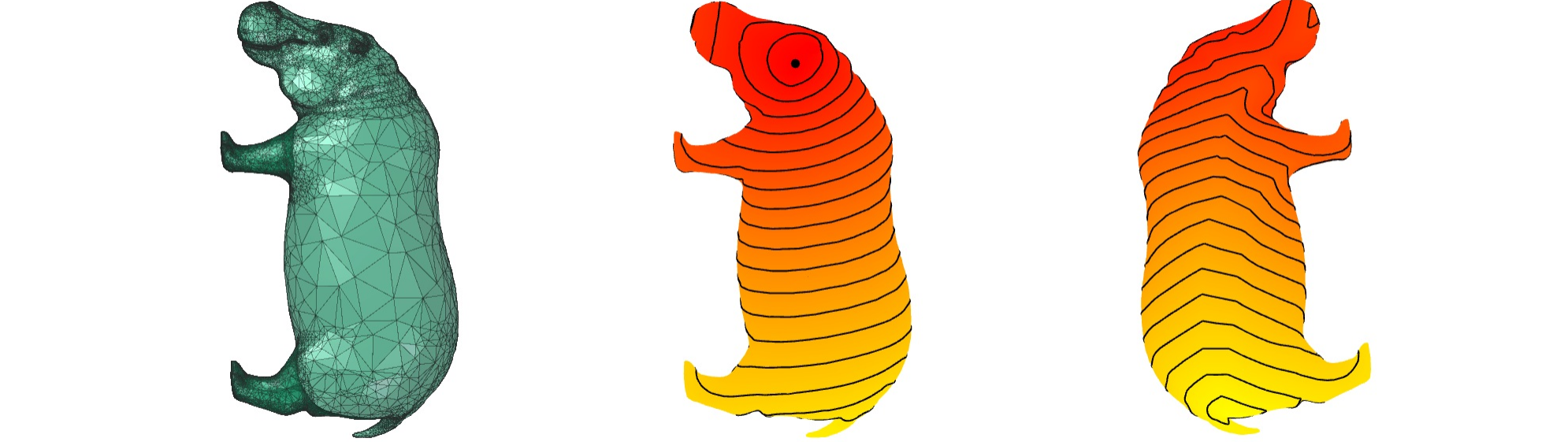} \\[10pt]
\includegraphics[width=0.9\textwidth]{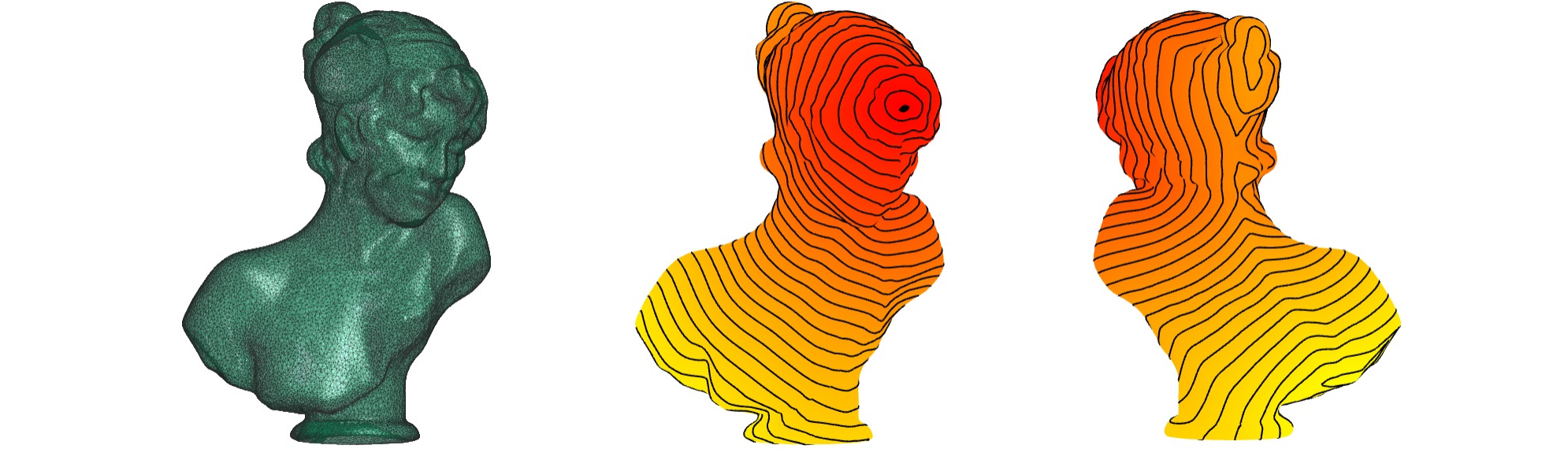}
\caption{CPHM solution of the Eikonal equation $\|\nabla_\S \phi\| = 1$ on various surfaces (in vertical order): Bunny, Bumpy Sphere, Pig, Hippo, and Sappho's Head. The numerical resolution for each case can be found in Table~\ref{tab:gallery_resolution}.}
\label{fig:gallery}
\end{figure}

We present the Eikonal solution with a single source on various surfaces in Fig.~\ref{fig:gallery}: Bunny, Bumpy Sphere, Pig, Hippo, and Sappho's Head. The corresponding number of grid points (inside the computational band) is provided in Table~\ref{tab:gallery_resolution}. Each row shows a different surface geometry, with the left column displaying the surface mesh, and the middle and right columns illustrating the computed geodesic distance via isolines from different viewpoints. The smoothness and density of the contours demonstrate the effectiveness of the method in capturing intrinsic distances across a variety of topologies and surface complexities.

\begin{table}[htbp]
    \centering
    \begin{tabular}{l c}
        \toprule
        Example & Length($\phi_h$) \\
        \midrule
        Bunny & 108,643 \\
        Bumpy sphere & 128,591 \\
        Pig & 278,573 \\
        Hippo & 396,311 \\
        Sappho's head \cite{Thingi10K} & 316,553 \\
        \bottomrule
    \end{tabular}
    \caption{Summary of grid resolutions used in Fig.~\ref{fig:gallery}.}
    \label{tab:gallery_resolution}
\end{table}

\section{Concluding Remarks}
In this paper, we presented the Closest Point Heat Method (CPHM), a novel approach for solving the surface Eikonal equation on general smooth surfaces. By extending the heat method to an embedding framework using the closest point methodology, CPHM overcomes several limitations of traditional techniques that rely on surface meshes or parametrizations. Our method enables intrinsic geodesic distance computations without requiring explicit surface discretizations, making it particularly effective for implicit surfaces or data represented as level sets or point clouds.

The key strengths of CPHM lie in its simplicity, mesh-free nature, and compatibility with standard finite difference tools on Cartesian grids. Through numerical experiments, we demonstrated the accuracy and convergence of the method on benchmark geometries and illustrated its applicability to complex shapes. The method maintains robustness even in the presence of geometric irregularities, while preserving the desirable properties of the original heat method, such as efficiency and stability.

Several promising directions remain for future research. A primary avenue is the enhancement of numerical accuracy and efficiency through adaptive or higher-order discretizations, as well as through the use of CPHM to generate high-quality initializations, for example, for learning-based methods \cite{smith2020eikonet,bin2021pinneik}. These improvements can be particularly beneficial in geometrically complex regions or learning frameworks where precise intrinsic information is crucial. Another central challenge is the extension of CPHM to anisotropic settings, where the speed function depends not only on position but also on direction. The anisotropic Eikonal equation arises in applications such as image processing, medical imaging, and fiber tractography, where wavefront propagation must align with directional features of the medium. Formally, it takes the form $\|A(\bm{y}) \nabla_{\S} \phi(\bm{y})\| = 1$, where $A(\bm{y})$ is a position-dependent tensor encoding local anisotropy. At present, our method is limited to the isotropic case with constant unit speed ($F(\bm{y}) = 1$), and even the extension to general spatially varying $F(\bm{y})$ remains an open challenge. Addressing the anisotropic case would require fundamental modifications to the heat flow approximation and projection operations, making it a natural but ambitious direction for future work.

\begin{appendix}

\section{Implementation of inhomogeneous Neumann boundary condition}\label{sec:NBC}

When applying the CPHM to an open surface $\S$ in Sec.~\ref{sec:open_surfaces}, we impose inhomogeneous Neumann BCs for the surface heat and Poisson equations. To the best of our knowledge, extending the CPM to handle general boundary conditions (besides homogeneous Neumann and Dirichlet BCs) remains an open challenge. Here, we generalize the modified closest point function (for open surfaces) proposed in \cite{macdonald2011solving} to treat an inhomogeneous Neumann BC. First, we first introduce the terminologies that lead to the modified closest point function.

\begin{figure}[t]
\centering
\tdplotsetmaincoords{60}{110}

\pgfmathsetmacro{\radius}{1}
\pgfmathsetmacro{\thetavec}{0}
\pgfmathsetmacro{\phivec}{0}

\begin{tikzpicture}[scale=3,tdplot_main_coords]
\draw[ultra thick,->] (1,0,0) -- (1.5,0,0) node[anchor=north]{$\textbf{N}$};
\draw[ultra thick,->] (1,0,0) -- (1,0,-0.4) node[anchor=north]{$\textbf{n}$};
\draw[ultra thick,->] (1,0,0) -- (1,0.4,0) node[anchor=north west]{$\textbf{T}$};
\tdplotsetthetaplanecoords{\phivec}

\draw[thick, dashed,tdplot_rotated_coords] (\radius,0,0) arc (0:90:\radius) ;

\draw[thick] (\radius,0,0) arc (0:360:\radius) node[pos=0.2,anchor=north west]{$\partial\mathcal{S}$};
\shade[ball color=blue!10!white,opacity=0.4] (1cm,0) arc (0:-180:1cm and 5mm) arc (180:0:1cm and 1cm) node[pos=0.3, anchor=north west, opacity=1.0]{$\mathcal{S}$};
\end{tikzpicture}\hfill~\begin{tikzpicture}[scale=3]
  \draw[thick, dashed] (1,0) arc[start angle=0, end angle=90, radius=1] node[anchor=south west]{$\mathcal{S}$};
  \draw[ultra thick,->] (1,0) -- (1.5,0) node[anchor=north]{$\textbf{N}$};
  \draw[ultra thick,->] (1,0) -- (1,-0.5) node[anchor=north]{$\textbf{n}$};

  \draw (1.0,  0.0) node[circle, fill=black, inner sep=.04cm]{}node[anchor=east]{$\textrm{cp}\left(\bm{x}_g\right)$};;
  \draw (0.7, -0.3) node[circle, fill=black, inner sep=.04cm]{}node[anchor=south east]{$\bm{x}_g$};
  \draw[dashed,-] (0.7,-0.3)--(1.0,0.0);
  \draw (1.3, 0.3,0.0) node [circle, fill=black, inner sep=.04cm]{}node[anchor=north west]{$2\textrm{cp}\left(\bm{x}_g\right)-\bm{x}_g$};
  \draw[dashed,-] (1.3,0.3)--(1.0,0.0);
  \draw (0.97,0.22) node [circle, fill=black, inner sep=.04cm]{}node[anchor=south east]{$\bar{\textrm{cp}}\left(\bm{x}_g\right)$};;
  \draw[dashed,-] (0.97,0.22)--(1.3,0.3);
  \foreach \i in {-90,-89,...,89}
  {
    \ifthenelse{\i>-90}{
        \draw[dashed]
          ({cos(\i)},{0.5*sin(\i)}) -- ({cos(\i+1)},{0.5*sin(\i+1)});
    }{
        \draw[dashed]
          ({cos(\i)},{0.5*sin(\i)}) -- ({cos(\i+1)},{0.5*sin(\i+1)}) node[anchor=north west]{$\partial\mathcal{S}$}; 
    }
    
  } 
  
\end{tikzpicture}
\caption{Illustration of mirror point construction near the boundary $\partial\S$ of the hemisphere $\S$. (Left) The surface normal $\bm{n}$, boundary normal $\bm{N}$, and tangential direction $\bm{T}$ at a boundary point. (Right) A grid point $\bm{x}_g$, its closest point $\cp\left(\bm{x}_g\right)$, and the corresponding mirror point $2\cp\left(\bm{x}_g\right)-\bm{x}_g$. }\label{fig:ghost_point}
\end{figure}
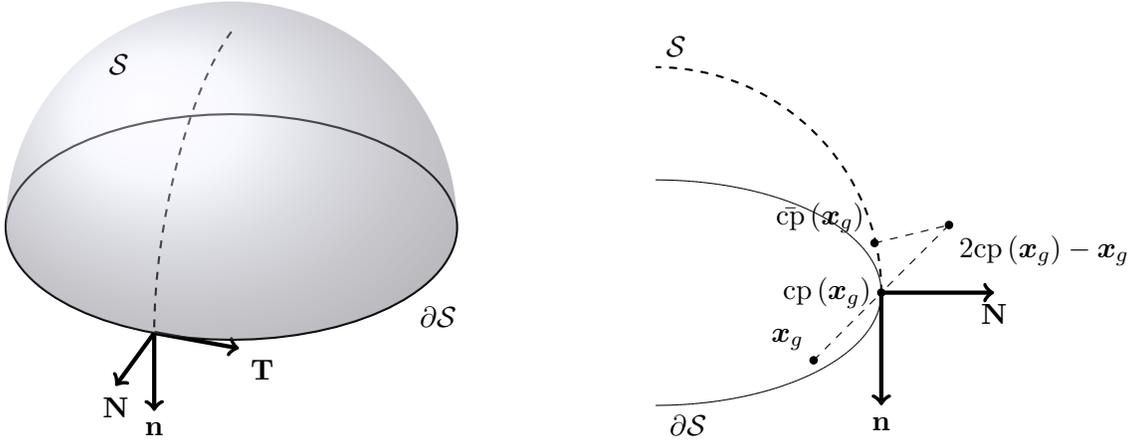

Let's consider a tubular neighbourhood $\mathcal{T}$ of an open surface $\S$ such that all grid points inside $\mathcal{T}$ have a unique closest point on $\S$. Let $v := u\circ \cp$ be the closest point extension of a surface function $u$ defined on $\S$. For any grid point $\bm{x}_g \in \mathcal{T}$, we have $v(\bm{x}_g) = v\left(\cp(\bm{x}_g)\right)$. Therefore, the closest point extension propagates the boundary values into the $\mathcal{T}$ along the normal directions to the boundary. In other words, when applied to an open surface, the closest point extension effectively imposes the homogeneous Neumann BC $\partial_n u = 0$ on $\partial S$. See Fig.~\ref{fig:ghost_point}.

Now, let's consider the modified closest point function from \cite{macdonald2011solving}
\begin{equation*}
    \cpbar\left({\bm{x}_g}\right) := \cp\left(\bm{x}_g + \cp(\bm{x}_g) - \bm{x}_g\right) = \cp\left(2\cp(\bm{x}_g) - \bm{x}_g \right) \,.
\end{equation*}
In here, $2\cp(\bm{x}_g) - \bm{x}_g$ is the ``mirror point" of $\bm{x}_g$ in the direction $\bm{x}_g - \cp(\bm{x}_g)$ across $\S$. Suppose that $\bm{x}_g - \cp(\bm{x}_g)$ is orthogonal to $\S$. Then both $\bm{x}_g$ and its mirror point $2\cp(\bm{x}_g) - \bm{x}_g$ share the same orthogonal projection, $\cp(\bm{x}_g)$, on $\S$. That is, $\cp(\bm{x}_g) = \cpbar(\bm{x}_g)$, which is true when $\bm{x}_g$ is away from the boundary $\partial S$. This naturally leads to a classification of boundary points in the embedding space. Specifically, we say that $\bm{x}_g$ is a boundary point if $\cp(\bm{x}_g) \in \partial \S$. Consequently, the vector $\bm{x}_g - \cp(\bm{x}_g)$ is not orthogonal to $\S$ and $\cp(\bm{x}_g) \neq \cpbar(\bm{x}_g)$. See Fig.~\ref{fig:ghost_point}. In this case, we identify $\bm{x}_g$ as a ghost point, whose value is treated as a boundary value. By replacing all instances of $\cp(\bm{x}_g)$ with $\cpbar(\bm{x}_g)$ in the discretization stencil, the modified closest point extension yields a second-order accurate discretization of the homogeneous Neumann BC. We let the modified interpolation matrix of order $q$ be $\bar{\mathsf{E}}_q$. Numerically, the modified closest point extension is enforced by the discrete constraint,
\begin{equation}\label{eqn:extrap_num}
    \mathsf{u} = \bar{\mathsf{E}}_q \mathsf{u} \,.
\end{equation}




Now, we utilize the modified closest point extension to handle the inhomogeneous Neumann BC $\partial_n u = g$ on $\partial S$. At a boundary point $\bm{x}_g$, let $\bm{T}$ and $\bm{N}$ be a unit tangent and normal vector of $\S$ such that $\bm{n} := \bm{T} \times \bm{N}$ is the unit outward normal of $\partial S$. Note that $\{\bm{T}, \bm{N}, \bm{n}\}$ is an orthonormal basis of $\mathbb{R}^3$. Let $\bm{w} = \frac{\bm{x}_g - \cp(\bm{x}_g)}{||\bm{x}_g - \cp(\bm{x}_g)||}$. Since $\bm{w}$ is orthogonal to $\bm{T}$, we have
\begin{equation}\label{eqn:w_decomp}
\bm{w} = \langle \bm{w}, \bm{n} \rangle \, \bm{n} + \langle \bm{w}, \bm{N} \rangle \bm{N} \,.
\end{equation}
From \eqref{eqn:w_decomp}, we may decompose the directional derivative in the $\bm{w}$-direction as
$$\partial_w = \langle \bm{w}, \bm{n} \rangle \, \partial_n + \langle \bm{w}, \bm{N} \rangle \, \partial_N \,.$$ 
Apply this to $u\big(\cp(\bm{x}_g))$, we have
\begin{equation}\label{eqn:directional_deriv}
\partial_w \, u\big( \cp(\bm{x}_g) \big) = \langle \bm{w}, \bm{n} \rangle \, g\big(\cp(\bm{x}_g)\big) \,.
\end{equation}
where we use the given BC $\partial_n u = g$ at $\bm{x} = \cp(\bm{x}_g)$, and the fact that $u\big(\cp(\bm{x})\big)$ is constant along the normal direction of $S$, hence $\partial_N u\big(\cp(\bm{x}_g)\big) = 0$. Next, we approximate the LHS of \eqref{eqn:directional_deriv} by the central difference,
\begin{equation}
\partial_w u\big(\cp(\bm{x}_g)\big) \approx \frac{u(\bm{x}_g) - u\big(\cpbar(\bm{x}_g)\big)}{2||\bm{x}_g - \cp(\bm{x}_g)||} \,.
\end{equation}
This yields
\[
\frac{u(\bm{x}_g) - u(\bar{\cp}(\bm{x}_g))}{2\|\bm{x}_g - \bar{\cp}(\bm{x}_g)\|} \approx  \langle \bm{w}, \bm{n} \rangle \, g(\cp(\bm{x}_g)),
\]
which eventually leads to a second-order accurate extrapolation formula:
\begin{equation}\label{eqn:NBC_extrapolation}
u(\bm{x}_g) = u\left(\bar{\cp}(\bm{x}_g)\right) + \langle \bm{x} - \cp(\bm{x}), \, \bm{n} \rangle\, g(\cp(\bm{x})).
\end{equation}
This modifies the constraint \eqref{eqn:extrap_num} to
\begin{equation}\label{eqn:extrap_num2}
\mathsf{u} = \bar{\mathsf{E}}_q \mathsf{u} + \mathsf{g} \,,
\end{equation}
where $\mathsf{g}$ is the vector that encodes the inhomogeneous Neumann BC (the second term of the RHS in \eqref{eqn:NBC_extrapolation}) at the boundary grid points. To demonstrate the second-order convergence with the mirror point treatment of the inhomogeneous Neumann BC, we consider a shifted Poisson equation \eqref{eqn:surf_spoisson}  on the upper hemisphere $\S$,
\begin{equation*}
    \Delta_\S u - u = f \enspace \mbox{in} \enspace \S \,, \qquad \partial_n u = g \enspace \mbox{on} \enspace \partial \S \,,
\end{equation*}
In spherical coordinate $(\varphi, \theta) \in [0,2\pi)\times(0, \pi)$, the exact solution is chosen to be $u(\varphi, \theta) = \sin\theta \cos\theta \cos\varphi$. Then we have $f=-7u$ and $g = -\cos\varphi$. We present the relative error and the convergence order in Table~\ref{tab:NBC_convergence}.

\begin{table}[htbp]
  \centering
  \begin{tabular}{ccc}
    \toprule
    $\Delta x$ &
    $\displaystyle\frac{\|u - u_h\|_\infty}{\|u\|_\infty}$ &
    Order \\
    \midrule
    $0.1$ & 6.6396E-3 & ---  \\
    $0.05$ & 1.8217E-3 & 1.8658 \\
    $0.025$ & 4.7954E-4 & 1.9256 \\
    $0.0125$ & 1.2362E-4 & 1.9557 \\
    \bottomrule
  \end{tabular}
  \caption{Second-convergence of the inhomogeneous Neumann problem}
  \label{tab:NBC_convergence}
\end{table}

\end{appendix}

\section*{Acknowledgement}
The research of Byungjoon Lee was supported by the Catholic University of Korea, Research Fund. The authors are grateful to anonymous reviewers for their careful reading and valuable comments. 

\bibliographystyle{vancouver}
\bibliography{ref}

\end{document}